\documentclass{article}

\usepackage{arxiv}
\usepackage[table,xcdraw]{xcolor}
\usepackage{amsmath,amsthm,amssymb,amscd}
\usepackage[all,cmtip]{xy}
\usepackage[utf8]{inputenc} 
\usepackage[T1]{fontenc}    
\usepackage{hyperref}       
\usepackage{url}            
\usepackage{booktabs}       
\usepackage{amsfonts}       
\usepackage{nicefrac}       
\usepackage{microtype}      
\usepackage{lipsum}
\usepackage{mathrsfs}
\usepackage{fancyhdr}
\usepackage{color}
\usepackage{subfig}
\usepackage{graphicx}
\usepackage{textcomp}
\usepackage{ytableau}
\usepackage{tikz}
\usetikzlibrary{calc,intersections,through,backgrounds,patterns}
\newtheorem{exam.}{Example}
\newtheorem{def.}{Definition}

\newtheorem{prop.}{Proposition}

\newtheorem{theorem}{Theorem}
\newtheorem*{remark}{Remark}

\usepackage{algorithm}
\usepackage{algpseudocode}

\newcommand{\id}{{\rm{id}}}

\newcommand{\Hom}{{\rm{Hom}}}

\newcommand{\bfx}{\mathbf{x}}

\title{Cellular Sheaves on Higher-Dimensional Structures}
\author{
Chuan-Shen Hu \\ 
Division of Mathematical Sciences\\
School of Physical and Mathematical Sciences\\
Nanyang Technological University \\
Singapore 637371 \\
\texttt{chuanshenhu.official@gmail.com} \\
}

\begin{document}
\maketitle

\begin{abstract} 
Defining cellular sheaves beyond graph structures, such as on simplicial complexes containing higher-dimensional simplices, is an essential and intriguing topic in topological data analysis (TDA) and the development of sheaf neural networks. In this paper, we explore methods for constructing non-trivial cellular sheaves on spaces that include structures of dimension greater than one. This extends the focus from 0- or 1-dimensional components, such as vertices and edges, to elements like triangles, tetrahedra, and other higher-dimensional simplices within a simplicial complex. We develop a unified framework that incorporates both geometric and algebraic approaches to modeling such complex systems using cellular sheaf theory. Motivated by the geometric and physical insights from anisotropic network models (ANM), we first introduce constructions that define sheaf structures whose 0-th sheaf Laplacians recover classical ANM Hessian matrices. The higher-dimensional sheaf Laplacians in this setting encode additional patterns of multi-way interactions. In parallel, we propose an algebraic framework based on commutative algebra and ringed spaces, where sheaves of ideals and modules are used to define sheaf structures in a combinatorial and algebraically grounded manner. These two perspectives—the geometric-physical and the algebraic—offer complementary strengths and together provide a versatile framework for encoding structural relationships and analyzing multi-scale data over simplicial complexes.
\end{abstract}
\section{Introduction}

Cellular sheaves provide a unifying framework for encoding and analyzing local-to-global relationships across structured domains such as graphs, simplicial complexes, and cell complexes. In recent years, cellular sheaf theory has found growing applications in topological data analysis (TDA)~\cite{kashiwara2018persistent,ghrist2022network,hansen2021opinion,curry2015topological,curry2016discrete,XiaoqiWei2024FODS}, signal processing~\cite{robinson2014topological,hansRiess2022,riess2022lattice}, sensor networks~\cite{robinson2014topological},  mechanical rigidity theory~\cite{cooperband2024cellular,cooperband2023cosheaf,cooperband2024equivariant,cooperband2023towards}, and the design of sheaf-based neural networks~\cite{hansen2020sheaf,he2023sheaf,braithwaite2024heterogeneous,duta2024sheaf,ayzenberg2025sheaf,barbero2022sheaf_attention,caralt2024joint,battiloro2023tangent,bodnar2022neural,bodnar2021weisfeiler,papamarkou2024position,hajij2025copresheaftopologicalneuralnetworks}. However, from an implementation and application perspective, much of the existing literature focuses primarily on sheaves defined over 0- or 1-dimensional structures, such as graphs or simplicial skeletons consisting only of vertices and edges. Extending these constructions to include higher-dimensional simplices, such as triangles, tetrahedra, and other complex structures, is a natural and necessary progression for modeling increasingly intricate data and interactions in geometry, physics, and data science. 

Although approaches such as labeled simplicial complexes, based on physical considerations in molecular systems, have been proposed to extend sheaf structures to higher-dimensional relationships within simplicial complex representations of atomic configurations~\cite{XiaoqiWei2024FODS,hayes2025persistent}, the discovery of diverse construction methods for cellular sheaves on simplicial complexes with higher-dimensional simplices remains an intriguing and largely unexplored area. This includes sheaves that take values in vector spaces or rings, which may offer alternative insights into complex data structures and interactions.

In this work, we explore methods for constructing non-trivial cellular sheaves on simplicial complexes of dimension greater than one. The goal is to develop a systematic approach for extending sheaf structures beyond graphs to encode higher-order relationships. This line of inquiry is particularly relevant to topological learning and geometric machine learning, where higher-dimensional interactions are prevalent but often challenging to capture.

We propose and investigate three complementary approaches for defining sheaf structures on simplicial complexes. The first approach is inspired by anisotropic network models (ANM), a class of Elastic Network Models (ENM) originally developed and widely applied in structural biology for coarse-grained modeling of molecular fluctuations~\cite{Atilgan:2001,doruker2000dynamics,eyal2006anisotropic,Opron:2014,Xia:2015multiscale,xia2014identifying}, as well as their recent formulation in terms of cellular sheaves~\cite{hu2024physicsinformedsheafmodel}. Building on this foundation, we explore the possibility of defining stalks and restriction maps not only on edges but also on 2-dimensional simplices and their associated edge-face relations. In particular, by interpreting ANM Hessians as 0-th sheaf Laplacians and extending the framework in~\cite{hu2024physicsinformedsheafmodel}, we construct a cellular sheaf on 2-dimensional simplicial complexes that encodes geometric relationships among vertices, edges, and faces. This sheaf preserves the original ANM Hessian structure while revealing higher-dimensional edge interactions through its 1-st sheaf Laplacian.

The second approach constructs higher-dimensional sheaves through the tensor product of sheaves, enabling the integration of geometric and physical information across different levels of the simplicial hierarchy. This tensor-based construction captures complex transition relationships between adjacent simplices and provides a flexible mechanism for enriching the overall sheaf structure.

The third approach shifts from geometric intuition to an algebraic perspective, drawing on tools from commutative algebra. Inspired by the theory of ringed spaces and sheaves of modules over posets, we define cellular sheaves of rings and modules using algebraic data such as ideals and localizations in polynomial rings. These constructions are particularly well-suited for exploring the combinatorial and algebraic geometry inherent in simplicial complexes. In particular, this approach is closely related to the theory of weighted homology on simplicial complexes \cite{dawson1990homology,bell2019weighted,ren2018weighted,meng2019weighted}, offering a sheaf-theoretic perspective that provides an alternative interpretation of the algebraic structure underlying weighted homology.

Taken together, these three paradigms—geometric-physical and algebraic—offer distinct yet complementary perspectives on the design and interpretation of cellular sheaves. They illustrate how topological and algebraic frameworks can be integrated to model multiscale interactions in both applied and theoretical contexts. The methods developed here provide a foundation for future work in sheaf-theoretic modeling and its applications in data analysis, mathematical physics, and related fields.

\section{Notations and Terminologies}

Typical notations for number systems are used throughout this paper, such as \( \mathbb{N} \subseteq \mathbb{Z} \subseteq \mathbb{Q} \subseteq \mathbb{R} \subseteq \mathbb{C} \), denoting the sets of positive integers, integers, rational numbers, real numbers, and complex numbers, respectively. The sets of all non-negative integers and non-negative real numbers are denoted by \( \mathbb{Z}_{\geq 0} \) (or \( \mathbb{N}_0 \)) and \( \mathbb{R}_{\geq 0} \), respectively. An arbitrary field (e.g., \( \mathbb{Q} \), \( \mathbb{R} \), \( \mathbb{C} \), etc.) is often denoted by \( \mathbb{F} \).

The underlying abstract simplicial complex—hereafter referred to simply as a \textit{simplicial complex}—for the proposed cellular sheaf models is assumed to be a finite abstract simplicial complex \( (K, \leq) \), where \( \leq \) denotes the face relation. Based on this face relation \( \leq \), the pair \( (K, \leq) \) forms a partially ordered set, or \textit{poset} for short. The vertex set of \( K \) is typically denoted by \( V = \{ 1, 2, \ldots, n \} \), where \( K \) consists of all non-empty subsets of \( V \) satisfying the property that if \( \sigma \in K \) and \( \tau \subseteq \sigma \), then \( \tau \in K \). The collection of all \( q \)-simplices with \( q \geq 0 \) (i.e., simplices of cardinality \( q+1 \)) is denoted by \( K_{(q)} \). Furthermore, we conventionally require that every vertex in \( V \) forms a \( 0 \)-simplex in \( K \); that is, \( V = K_{(0)} \). To simplify notation, we will often use \( i \in V \) to denote the \( 0 \)-simplex \( [i] \) and write \( i \leq \sigma \) for \( \sigma \in K \) if \( i \) corresponds to a vertex of \( \sigma \) (i.e., $i \in \sigma$).

With the natural total order on the vertices \( 1 < 2 < \cdots < n \), every \( q \)-simplex with \( q \in \mathbb{Z}_{\geq 0} \) is uniquely represented as an ordered tuple \( [i_0, \ldots, i_q] \) with \( i_0 < i_1 < \cdots < i_q \). Under this total order, the standard incidence sign function \( [\sigma : \tau] \) can be defined, inheriting the orientation of simplices as in homology theory. Specifically, \( [\sigma : \tau] = (-1)^k \) if \( \tau = [i_0, \ldots, i_q] \) and \( \sigma = [i_0, \ldots, \widehat{i_k}, \ldots, i_q] \), and \( [\sigma : \tau] = 0 \) otherwise.

In some circumstances, in addition to the combinatorial relations between faces, the underlying simplicial complex is also considered as a geometric object embedded in \( \mathbb{R}^d \). Specifically, each vertex $i \in V \subseteq K$ is accompanied by a $d$-dimensional real-valued vector $\mathbf{r}_i \in \mathbb{R}^d$. As a convention, the notation for a vector $\mathbf{r}$ in $\mathbb{R}^d$ represents a column vector, i.e., a $3 \times 1$ matrix. Examples in this paper frequently illustrate the case where \( d = 3 \), representing the 3D coordinates of the vertices.

In this paper, several notions and conventions from commutative algebra are also adopted. Specifically, every \textit{ring} refers to a commutative ring with identity, and every \textit{ring homomorphism} is assumed to preserve the identity elements of the respective rings. Given a subset \( S \subseteq R \), the ideal generated by \( S \) is denoted by \( (S) \). In the case where \( S = \{ s_1, \dots, s_k \} \) is finite, we write \( (s_1, \dots, s_k) \). The \textit{prime spectrum} of a ring \( R \), denoted by \( \operatorname{Spec}(R) \), is the set of all prime ideals of \( R \). This set is often simply referred to as the \textit{spectrum} of \( R \). The polynomial ring in \( n \) indeterminates over a field \( \mathbb{F} \) is denoted by \( \mathbb{S}_n = \mathbb{F}[x_1, \ldots, x_n] \), or simply by \( \mathbb{S} \) when the context makes the choice of indeterminates clear.

\section{Cellular Sheaves on Simplicial Complexes}
\label{Section: Cellular Sheaves on Simplicial Complexes}

A basic familiarity with category theory is assumed throughout this paper. Specifically, a \textit{category} \( \textsf{C} \) consists of a class of \textit{objects}, denoted by \( \mathrm{Obj}(\textsf{C}) \), and a class of \textit{morphisms} (also called \textit{arrows}) between these objects. For any two objects \( A, B \in \mathrm{Obj}(\textsf{C}) \), the set of morphisms from \( A \) to \( B \) is denoted by \( \Hom_{\textsf{C}}(A, B) \). A morphism \( f \in \Hom_{\textsf{C}}(A, B) \) is typically written as \( f: A \rightarrow B \). The \textit{opposite category} of a category \( \textsf{C} \), denoted by \( \textsf{C}^{\mathrm{op}} \), consists of the same objects as in \( \textsf{C} \), but with all morphisms reversed; that is, \( \Hom_{\textsf{C}^{\mathrm{op}}}(A, B) := \Hom_{\textsf{C}}(B, A) \) for any objects \( A \) and \( B \).

Typical categories considered in this paper include the category of topological spaces \( \textsf{Top} \), the category of Abelian groups \( \textsf{Ab} \), the category of rings \( \textsf{Ring} \), the category of \( R \)-modules \( \textsf{Mod}_R \) over a fixed ring \( R \), and the category of \( F \)-vector spaces \( \textsf{Vect}_F \) over a fixed field \( \mathbb{F} \), among others. Furthermore, from a computational perspective, the category \( \textsf{vect}_{\mathbb{F}} \) of finite-dimensional \( \mathbb{F} \)-vector spaces is primarily considered, rather than the category \( \textsf{Vect}_{\mathbb{F}} \).

\paragraph{Formal definitions}
For a given ring \( R \), a \textit{cellular sheaf} of \( R \)-modules over a simplicial complex \( (K, \leq) \) is defined as a (covariant) functor \( \mathcal{F} : (K, \leq) \rightarrow \textup{\textsf{Mod}}_R \). This functor assigns to each simplex \( \sigma \in K \) an \( R \)-module \( \mathcal{F}_\sigma \), and to each face relation \( \sigma \leq \tau \), an \( R \)-module homomorphism \( \mathcal{F}_{\sigma,\tau} : \mathcal{F}_\sigma \rightarrow \mathcal{F}_\tau \), such that for any \( \sigma \leq \tau \leq \eta \), the composition satisfies \( \mathcal{F}_{\tau,\eta} \circ \mathcal{F}_{\sigma,\tau} = \mathcal{F}_{\sigma,\eta} \). The \( R \)-module \( \mathcal{F}_\sigma \) is called the \textit{stalk} at \( \sigma \), and the map \( \mathcal{F}_{\sigma,\tau} \) is called the \textit{restriction map} from \( \mathcal{F}_\sigma \) to \( \mathcal{F}_\tau \).

\paragraph{Sheaf cohomology} 
The \textit{sheaf cohomology} and \textit{sheaf Laplacian} associated with a cellular sheaf \( \mathcal{F} : (K, \leq) \rightarrow \textup{\textsf{Mod}}_R \) are also studied in this paper. Specifically, for every non-negative integer $q \in \mathbb{Z}_{\geq 0}$, the $q$-\textit{th cochain space} of $\mathcal{F}$ is defined as the direct sum of the stalks of $\mathcal{F}$ on $q$-simplices. That is,
\begin{equation*}
C^q(K;\mathcal{F}) = \prod_{\sigma \in K_{(q)}} \mathcal{F}_\sigma.
\end{equation*}
\begin{remark}
Since \( K \) is assumed to be finite, the direct product \( \prod_{\sigma \in K_{(q)}} \mathcal{F}_\sigma \) and direct sum \( \bigoplus_{\sigma \in K_{(q)}} \mathcal{F}_\sigma \) coincide.   
\end{remark}
Elements in $C^q(K;\mathcal{F})$ are usually expressed by tuples $(\mathbf{x}_\sigma)_{\sigma \in K_q}$ with $\mathbf{x}_\sigma \in \mathcal{F}_\sigma$ for every $\sigma \in K_q$. The $q$-\textit{th coboundary map} $\delta^q: C^q(K;\mathcal{F}) \rightarrow C^{q+1}(K;\mathcal{F})$ is thus defined as the $R$-module homomorphism extended by the following assignment:
\begin{equation*}
\delta^q|_{\mathcal{F}_\sigma} = \sum_{\tau \in K_{(q+1)}} [\sigma:\tau] \cdot \mathcal{F}_{\sigma, \tau}.     
\end{equation*}
A standard argument in (co)homology theory shows that, for every \( q \geq 0 \), the composition \( \delta^{q+1} \circ \delta^q \) of the coboundary maps from \( C^q(K, \mathcal{F}) \) to \( C^{q+2}(K, \mathcal{F}) \) satisfies \( \delta^{q+1} \circ \delta^q = 0 \). Consequently, the \( q \)-th sheaf cohomology is defined as the quotient $R$-module
\begin{equation*}
H^q(K, \mathcal{F}) = \frac{\ker(\delta^q)}{\operatorname{im}(\delta^{q-1})}.
\end{equation*}
The $0$-th sheaf cohomology of a cellular sheaf \( \mathcal{F} : (K, \leq) \rightarrow \textup{\textsf{Mod}}_R \) is referred to as the \textit{space of global sections} of \( \mathcal{F} \), and is typically denoted by \( \Gamma(K, \mathcal{F}) \). In sheaf theory, this space consists of \textit{harmonic signals} on the simplices that can be coherently assembled into a global signal~\cite{robinson2014topological, hansen2019toward, bodnar2022neural}. In our context, the space of global sections is mathematically represented as the inverse limit
\begin{equation}\label{Eq. global section space as an inverse limit}
\Gamma(K;\mathcal{F}) = \varprojlim_{\sigma \in (K, \leq)} \mathcal{F}_\sigma = \left\{ \left. (\bfx_\sigma)_{\sigma \in K} \in \prod_{\sigma \in K} \mathcal{F}_\sigma \ \right| \ \mathcal{F}_{\sigma, \eta}(\bfx_\sigma) = \mathcal{F}_{\tau, \eta}(\bfx_\tau) \text{ if } \sigma \leq \eta \text{ and } \tau \leq \eta \right\}.    
\end{equation}
Because every \( 0 \)-simplex in the simplicial complex \( K \) is a minimal element in the poset \( (K, \leq) \), the space of global sections is canonically isomorphic, as an \( R \)-module, to
\begin{equation}\label{Eq. global section space as an inverse limit-simple form}
\Gamma(K, \mathcal{F}) \simeq \left\{ (v_\sigma)_{\sigma \in K_{(0)}} \in \prod_{\sigma \in K_{(0)}} \mathcal{F}_\sigma \ \middle| \ \mathcal{F}_{\sigma,\tau}(v_\sigma) = \mathcal{F}_{\sigma',\tau}(v_{\sigma'}) \text{ for all } \tau \in K_{(1)} \text{ such that } \sigma, \sigma' \leq \tau \right\}.    
\end{equation} 
Using a matrix representation, and assuming that \( K \) is finite, the coboundary maps of a cellular sheaf \( \mathcal{F} : (K, \leq) \rightarrow \textup{\textsf{Mod}}_R \) can be expressed as block matrices with entries in \( R \) of the form:
\begin{equation}
\label{Eq. Coboundary matrix}
\mathbf{C}^q = \begin{bmatrix}
[\sigma_1:\tau_1] \cdot \mathcal{F}_{\sigma_1, \tau_1} & [\sigma_2:\tau_1] \cdot \mathcal{F}_{\sigma_2, \tau_1} & \cdots & [\sigma_n:\tau_1] \cdot \mathcal{F}_{\sigma_n, \tau_1}\\
[\sigma_1:\tau_2] \cdot \mathcal{F}_{\sigma_1, \tau_2} & [\sigma_2:\tau_2] \cdot \mathcal{F}_{\sigma_2, \tau_2} & \cdots & [\sigma_n:\tau_2] \cdot \mathcal{F}_{\sigma_n, \tau_2}\\
\vdots & \vdots & \ddots & \vdots \\
[\sigma_1:\tau_m] \cdot \mathcal{F}_{\sigma_1, \tau_m} & [\sigma_2:\tau_m] \cdot \mathcal{F}_{\sigma_2, \tau_m} & \cdots & [\sigma_n:\tau_m] \cdot \mathcal{F}_{\sigma_n, \tau_m}\\
\end{bmatrix},
\end{equation}
where $K_{(q)} = \{ \sigma_1, \sigma_2, ..., \sigma_n \}$ and $K_{(q+1)} = \{ \tau_1, \tau_2, ..., \tau_m \}$. In particular, if \( \mathcal{F}_\sigma = \mathbb{F}^{d_q} \) for all \( \sigma \in K_{(q)} \) and \( \mathcal{F}_\tau = \mathbb{F}^{d_{q+1}} \) for all \( \tau \in K_{(q+1)} \), then the coboundary matrix \( \mathbf{C}^q \) is an \( (m d_{q+1}) \times (n d_q) \) matrix with entries in \( \mathbb{F} \). In the representation of the matrix~$\mathbf{C}^q$ in~\eqref{Eq. Coboundary matrix}, we adopt the convention that $\mathcal{F}_{\sigma,\tau}$ is defined as the zero map whenever $\sigma$ is not a face of $\tau$. This choice is consistent with the coefficient~$[\sigma:\tau]$.

Sheaf cohomology and the sheaf Laplacian are closely connected. Specifically, for any cellular sheaf of finite-dimensional \( \mathbb{R} \)-vector spaces on a finite simplicial complex---that is, a functor \( \mathcal{F} : (K, \leq) \rightarrow \textup{\textsf{vect}}_\mathbb{R} \)---the \( q \)-th \textit{combinatorial Hodge Laplacian} is defined on the cochain complex
\begin{equation*}
\xymatrix@+0.0em{
\cdots
\ar[r]
& C^{q-1}(K;\mathcal{F})
\ar[r]^{\delta^{q-1}}
& C^{q}(K;\mathcal{F})
\ar[r]^{\delta^{q}}
\ar@/^1pc/[l]^{(\delta^{q-1})^*}
& C^{q+1}(K;\mathcal{F})
\ar@/^1pc/[l]^{(\delta^{q})^*}
\ar[r]
& \cdots
}
\end{equation*}
as the linear operator $\Delta^q_{\mathcal{F}} = (\delta^{q})^* \circ \delta^q + \delta^{q-1} \circ (\delta^{q-1})^*$, which maps \( C^q(K; \mathcal{F}) \) to itself. Here, \( (\delta^\bullet)^* \) denotes the adjoint \( \mathbb{R} \)-linear transformation of \( \delta^\bullet \) with respect to the standard inner product on cochain spaces.

\begin{theorem}
Let $\mathcal{F}: (K,\leq) \rightarrow \textup{\textsf{vect}}_{\mathbb{R}}$ be a cellular sheaf on a simplicial complex $K$, and let $\Delta^q_\mathcal{F}$ be the $q$-th Hodge Laplacian of $\mathcal{F}$. Then, $H^q(K;\mathcal{F}) \simeq \ker(\Delta^q)$. In particular, 
\begin{equation}
\label{Eq. Global section space and the kernel of 0th Laplacian}
\Gamma(K;\mathcal{F}) = H^0(K;\mathcal{F}) = \ker(\Delta^0) = \ker((\delta^0)^* \cdot \delta^0).
\end{equation}
\end{theorem}
\begin{proof}
For a proof of a more general version of the theorem in the setting of cellular sheaves on cellular complexes, see~\cite[Theorem~3.1]{hansen2019toward}.
\end{proof}

In particular, if \( \mathcal{F}_\sigma = \mathbb{R}^{d_q} \) for all \( \sigma \in K_{(q)} \) and \( \mathcal{F}_\tau = \mathbb{R}^{d_{q+1}} \) for all \( \tau \in K_{(q+1)} \), then the \( q \)-th coboundary map \( \delta^q \) is represented by a coboundary matrix \( \mathbf{C}^q \), which is an \( (m d_{q+1}) \times (n d_q) \) matrix with entries in \( \mathbb{R} \), where \( n = |K_{(q)}| \) and \( m = |K_{(q+1)}| \). The adjoint of \( \delta^q \) corresponds to the transpose matrix \( (\mathbf{C}^q)^\top \). Consequently, the matrix representing the \( q \)-th sheaf Laplacian is given by
\begin{equation}\label{Eq. Sheaf Laplacian}
\mathbf{L}^q_{\mathcal{F}} = 
(\mathbf{C}^q)^\top \cdot \mathbf{C}^q =
\begin{bmatrix}
L_{11} & L_{12} & \cdots & L_{1n} \\
L_{21} & L_{22} & \cdots & L_{2n} \\
\vdots & \vdots & \ddots & \vdots \\
L_{n1} & L_{n2} & \cdots & L_{nn} \\
\end{bmatrix},
\end{equation}
where the matrix is block-indexed by the set \( K_{(q)} = \{ \sigma_1, \sigma_2, \ldots, \sigma_n \} \), consisting of \( n \) simplices. Each block \( L_{ij} \) is a \( d_q \times d_q \) matrix with entries in \( \mathbb{R} \).

Based on the framework of cellular sheaves, along with their associated coboundary maps and Laplacian matrices, classical mathematical objects---such as graph Laplacians (also known as Kirchhoff matrices)---can be recovered as special cases of the sheaf Laplacian.

\begin{exam.}[Graph Laplacians]\label{Example: Graph Laplacians}
Let \( G = (V, E) \) be a finite, unweighted, simple graph, where \( V = \{ 1, 2, ..., n \} \) denotes the set of vertices and \( E \subseteq \{ \{i, j\} \mid i, j \in V, i \neq j \} \) denotes the set of edges. Then the graph Laplacian of $G$ is an $n \times n$ matrix $\mathbf{L} = (L_{ij})$ that is defined by
\begin{equation*}
L_{ij} =
\begin{cases}
\deg(i) &\text{if} \quad i = j,  \\
-1 &\text{if} \quad i \neq j \text{ and } \{ i, j \} \in E,\\
0  &\text{if} \quad i \neq j \text{ and } \{ i, j \} \notin E,
\end{cases}    
\end{equation*}
where $\deg(i)$ denotes the degree of vertex $i$ based on the graph structure. Then, the graph Laplacian \( \mathbf{L} \) can be identified as a special case of the sheaf Laplacian under the following construction. Specifically, viewing \( G \) as the 1-dimensional simplicial complex \( V \cup E \), define a sheaf \( \mathcal{F} : (G, \leq) \rightarrow \textup{\textsf{vect}}_{\mathbb{R}} \) by setting \( \mathcal{F}_i = \mathbb{R} \) for each vertex \( i \in V \), \( \mathcal{F}_{\{i, j\}} = \mathbb{R} \) for each edge \( \{i, j\} \in E \), and defining the restriction maps \( \mathcal{F}_{i,\{i, j\}} = \mathcal{F}_{j,\{i, j\}} = \mathrm{id}_{\mathbb{R}} \). Under this construction, the classical graph Laplacian \( L \) coincides with the degree-zero sheaf Laplacian: \( \mathbf{L} = \mathbf{L}_{\mathcal{F}}^0 \).
\end{exam.}
\begin{proof}
Consider the standard ordering \( 1 < 2 < \cdots < n \) of the vertices of \( G \). Then, every edge of \( G \) can be uniquely expressed as \( [i, j] \) with \( i < j \). Suppose \( e_1, \ldots, e_m \) are all the edges of \( G \), then the corresponding degree-zero coboundary matrix \( \mathbf{C}^0 \) is given by:
\begin{equation*}
\mathbf{C}^0 = \begin{bmatrix}
[1: e_1] & [2: e_1] & \cdots & [n: e_1] \\
[ 1: e_2] & [ 2: e_2] & \cdots & [ n: e_2] \\
\vdots & \vdots & \ddots & \vdots \\
[ 1: e_m] & [ 2: e_m] & \cdots & [ n: e_m] \\
\end{bmatrix}.    
\end{equation*}
Then, the $0$-th sheaf Laplacian of $\mathcal{F}$ is the $n \times n$ matrix
\begin{equation*}
\mathbf{L}_{\mathcal{F}}^0 := (\mathbf{C}^0)^\top \cdot \mathbf{C}^0 = \begin{bmatrix}
l_{11} & l_{12} & \cdots & l_{1n} \\
l_{21} & l_{22} & \cdots & l_{2n} \\
\vdots & \vdots & \ddots & \vdots \\
l_{n1} & l_{n2} & \cdots & l_{nn} \\
\end{bmatrix},   
\end{equation*}
where $l_{ii} = \sum_{k = 1}^n [i:e_k] \cdot [i: e_k] = \sum_{k = 1}^n [i:e_k]^2 = \deg(i)$ for $i \in \{ 1, 2, ..., n \}$, and $l_{ij} = \sum_{k = 1}^n [i:e_k] \cdot [j: e_k]$ for $i \neq j$ in $\{ 1, 2, ..., n \}$. Based on the orientation \( 1 < 2 < \cdots < n \), $[i:e_k] \cdot [j: e_k] = -1$ if $e_k$ contains $i, j$ as its endpoints, and $[i:e_k] \cdot [j: e_k] = 0$ for otherwise. Because $G$ is a simple graph, there is at most one edge in $E$ that has endpoints $i$ and $j$. This shows that, for \( i \neq j \) in \( V \), the off-diagonal entry \( l_{ij} \) of the Laplacian matrix satisfies \( l_{ij} = -1 \) if \( \{ i, j \} \in E \), and \( l_{ij} = 0 \) otherwise. Therefore, the sheaf Laplacian \( \mathbf{L}_{\mathcal{F}}^0 \) recovers the standard graph Laplacian \( \mathbf{L} \); that is, $\mathbf{L} = \mathbf{L}_{\mathcal{F}}^0$.
\end{proof}

The sheaf \( \mathcal{F} : (G, \leq) \rightarrow \textup{\textsf{vect}}_{\mathbb{R}} \) defined in Example~\ref{Example: Graph Laplacians} is called the \textit{constant sheaf} on the graph \( G \), and is typically denoted by \( \underline{\mathbb{R}} \). A similar construction can be carried out with an arbitrary vector space $V \in \textup{\textsf{Vect}}_{\mathbb{R}}$, yielding the constant sheaf $\underline{V}$ on $G$, such as $\underline{\mathbb{R}^3}$. In the following example, we consider a slight generalization of the constant sheaf $\underline{\mathbb{R}^3}$ on the graph~$G$.

\begin{exam.}\label{Example: R3 constant-like sheaf on Graphs}
Let \( G = (V, E) \) be a finite, unweighted, simple graph, where \( V = \{ 1, 2, \ldots, n \} \) is the set of vertices and \( E \subseteq \{ \{i, j\} \mid i, j \in V,\, i \neq j \} \) is the set of edges. For every \( \lambda \in \mathbb{R} \setminus \{ 0 \} \), the following construction defines a cellular sheaf \( \mathcal{G} : (G, \leq) \rightarrow \textup{\textsf{vect}}_{\mathbb{R}} \): (a) \( \mathcal{G}_{i} = \mathbb{R}^3 \) for every \( i \in V \); (b) \( \mathcal{G}_{\{ i, j \}} = \mathbb{R}^3 \) for every \( \{ i, j \} \in E \); (c) \( \mathcal{G}_{i, \{ i, j \}} = \mathcal{G}_{j, \{ i, j \}} = \lambda \cdot I_3 \) for every \( \{ i, j \} \in E \), where \( I_3 \) denotes the \( 3 \times 3 \) identity matrix. In particular, the $0$-th sheaf Laplacian of $\mathcal{G}$ is
\begin{equation}\label{Eq. The first tensor form of the sheaf Laplacian}
\mathbf{L}_{\mathcal{G}}^0 = \mathbf{L}_{\mathcal{F}}^0 \otimes (\lambda^2 \cdot I_3) = \mathbf{L}_{\mathcal{F}}^0 \otimes \begin{bmatrix}
\lambda^2 & 0 & 0 \\
0 & \lambda^2 & 0 \\ 
0 & 0 & \lambda^2 \\
\end{bmatrix} = \lambda^2 \cdot \left( \mathbf{L}_{\mathcal{F}}^0 \otimes I_3 \right),   
\end{equation}
where the cellular sheaf $\mathcal{F} = \underline{\mathbb{R}}$ is the constant sheaf on the graph $G$, as defined in Example~\ref{Example: Graph Laplacians}, and $\otimes$ denotes the Kronecker product of matrices.
\end{exam.}
\begin{proof}
We present a detailed computation of the sheaf Laplacian matrix associated with the sheaf~$\mathcal{G}$. Suppose \( e_1, \ldots, e_m \) are all the edges of \( G \), then the corresponding degree-zero coboundary matrix \( \mathbf{C}_{\mathcal{G}}^0 \) is given by:
\begin{equation}\label{Eq. the coboundary matrix of GNM}
\begin{split}
\mathbf{C}_{\mathcal{G}}^0 &= \begin{bmatrix}
[1: e_1] \cdot \mathcal{G}_{1, e_1} & [2: e_1] \cdot \mathcal{G}_{2, e_1} & \cdots & [n: e_1] \cdot \mathcal{G}_{n, e_1} \\
[ 1: e_2] \cdot \mathcal{G}_{1, e_2} & [ 2: e_2] \cdot \mathcal{G}_{2, e_2} & \cdots & [ n: e_2] \cdot \mathcal{G}_{n, e_2} \\
\vdots & \vdots & \ddots & \vdots \\
[ 1: e_m] \cdot \mathcal{G}_{1, e_m} & [ 2: e_m] \cdot \mathcal{G}_{2, e_m} & \cdots & [ n: e_m] \cdot \mathcal{G}_{n, e_m}\\
\end{bmatrix},    
\end{split}
\end{equation}
where $[i:e_k] \cdot \mathcal{G}_{i, e_k} = \pm (\lambda \cdot I_3)$ if $i$ is an endpoint of edge $e_k$, and $[i:e_k] \cdot \mathcal{G}_{i, e_k}$ is the zero map if $i$ is not an endpoint of edge $e_k$. Therefore, $\mathbf{C}_{\mathcal{G}}^0 = \mathbf{C}_{\mathcal{F}}^0 \otimes (\lambda \cdot I_3)$, where $\mathbf{C}_{\mathcal{F}}^0$ is the $0$-th coboundary matrix of the cellular sheaf $\mathcal{F} = \underline{\mathbb{R}}$ on $G$ (Example~\ref{Example: Graph Laplacians}). Then, the $0$-th sheaf Laplacian of $\mathcal{G}$ is the matrix
\begin{equation*}
\begin{split}
(\mathbf{C}_{\mathcal{G}}^0)^\top \cdot \mathbf{C}_{\mathcal{G}}^0 &= (\mathbf{C}_{\mathcal{F}}^0 \otimes (\lambda \cdot I_3))^\top \cdot (\mathbf{C}_{\mathcal{F}}^0 \otimes (\lambda \cdot I_3)) = ((\mathbf{C}_{\mathcal{F}}^0)^\top \otimes (\lambda \cdot I_3)^\top) \cdot (\mathbf{C}_{\mathcal{F}}^0 \otimes (\lambda \cdot I_3)) \\
&= ((\mathbf{C}_{\mathcal{F}}^0)^\top \otimes (\lambda \cdot I_3)) \cdot (\mathbf{C}_{\mathcal{F}}^0 \otimes (\lambda \cdot I_3)) = ((\mathbf{C}_{\mathcal{F}}^0)^\top \cdot \mathbf{C}_{\mathcal{F}}^0) \otimes (\lambda^2 \cdot I_3) = \mathbf{L}_{\mathcal{F}}^0 \otimes (\lambda^2 \cdot I_3),
\end{split}   
\end{equation*}
and this proves Equation \eqref{Example: R3 constant-like sheaf on Graphs}. In particular, if $\lambda, \gamma \in \mathbb{R}_{> 0}$ with $\gamma = \lambda^2$, then $\mathbf{L}_{\mathcal{G}}^0 = \mathbf{L}_{\mathcal{F}}^0 \otimes (\gamma \cdot I_3)$.
\end{proof}

From the perspective of Elastic Network Models (ENMs), the sheaf Laplacian matrices \( \mathbf{L}_{\mathcal{F}}^0 \) and \( \mathbf{L}_{\mathcal{F}}^0 \otimes I_3 \), as shown in Example~\ref{Example: R3 constant-like sheaf on Graphs}, play a crucial role in the Gaussian Network Model (GNM) and the expanded Gaussian Network Model (eGNM). These matrices are used to express the GNM-based potential energy of the underlying macromolecule, modeled as an elastic mass-and-spring network~\cite{bahar1998vibrational,bahar1997direct,Opron:2014,park2013coarse,LWYang:2008}. 

Alternatively, from a different perspective on analyzing molecular fluctuations and dynamics, the Anisotropic Network Model (ANM) has also been widely applied~\cite{Atilgan:2001}. Inspired by the ANM, a sheaf structure can likewise be defined on the graph representing the atomic system in \( \mathbb{R}^3 \). In particular, as shown in the following example (see Example~\ref{Example: ANM sheaf on Graphs}), the \textit{Hessian matrix} in the ANM model, which encodes the system's potential energy, coincides with the sheaf Laplacian associated with this construction~\cite{hu2024physicsinformedsheafmodel}. For dual or equivalent constructions from the perspective of rigidity theory and mechanical frameworks, see~\cite{cooperband2024cellular,cooperband2023cosheaf,cooperband2024equivariant,cooperband2023towards}.

\begin{exam.}\label{Example: ANM sheaf on Graphs}
Let \( G = (V, E) \) be a finite, unweighted, simple graph, where \( V = \{ 1, 2, \ldots, n \} \) is the set of vertices and \( E \subseteq \{ \{i, j\} \mid i, j \in V,\, i \neq j \} \) is the set of edges. Each vertex \( i \) is assigned a column vector \( \mathbf{r}_i \in \mathbb{R}^3 \), with coordinates denoted by \( (x_i^\circ, y_i^\circ, z_i^\circ) \). Represent each edge in \( E \) as an ordered pair \( [i, j] \) with \( i < j \), and consider a sheaf \( \mathcal{F} : (G, \leq) \rightarrow \textup{\textsf{vect}}_{\mathbb{R}} \) defined as follows:
\begin{itemize}
\item[\rm (a)] $\mathcal{F}_i = \mathbb{R}^3$ for every vertex $i \in V$.
\item[\rm (b)] $\mathcal{F}_{[i,j]} = \mathbb{R}$ for every vertex $i, j \in V$ and edge $[i, j]$.
\item[\rm (c)] $\mathcal{F}_{i,[i,j]} = \mathcal{F}_{j,[i,j]} = \frac{\gamma^{1/2}}{d_{ij}^\circ} \cdot (\mathbf{r}_j - \mathbf{r}_i)^\top$ for every vertex $i, j \in V$ and edge $[i, j]$.  
\end{itemize}
In this representation, the row vector \( (\mathbf{r}_j - \mathbf{r}_i)^\top \) is regarded as a \( 1 \times 3 \) matrix, defining a linear transformation from \( \mathbb{R}^3 \) to \( \mathbb{R} \) via multiplication with column vectors in \( \mathbb{R}^3 \). The scalars \( \gamma \in \mathbb{R}_{> 0} \) and \( d_{ij}^\circ \) represent the \textbf{spring constant} and the \textbf{equilibrium distance} between the coordinates \( \mathbf{r}_i \) and \( \mathbf{r}_j \), respectively. Then, the ANM Hessian matrix, 
\begin{equation*}
\label{Eq. The big Hessian matrix}
\mathbf{H}_{\rm ANM} = \begin{bmatrix}
H_{11} & H_{12} & \cdots & H_{1n}\\
H_{21} & H_{22} & \cdots & H_{2n}\\
\vdots & \vdots & \ddots & \vdots \\
H_{n1} & H_{n2} & \cdots & H_{nn}\\
\end{bmatrix}   
\end{equation*}
with $H_{ii} := -\sum_{j \neq i} H_{ij}$ for $i = 1, 2, ..., n$ and
\begin{equation*}
H_{ij} = \frac{-\gamma}{(d_{ij}^\circ)^2} \cdot \begin{bmatrix} 
(x_j^\circ-x_i^\circ)(x_j^\circ-x_i^\circ) & (x_j^\circ-x_i^\circ)(y_j^\circ-y_i^\circ) & (x_j^\circ-x_i^\circ)(z_j^\circ-z_i^\circ) \\ (x_j^\circ-x_i^\circ)(y_j^\circ-y_i^\circ) & (y_j^\circ-y_i^\circ)(y_j^\circ-y_i^\circ) & (y_j^\circ-y_i^\circ)(z_j^\circ-z_i^\circ) \\ (x_j^\circ-x_i^\circ)(z_j^\circ-z_i^\circ) & (y_j^\circ-y_i^\circ)(z_j^\circ-z_i^\circ) & (z_j^\circ-z_i^\circ)(z_j^\circ-z_i^\circ)  \\
\end{bmatrix}
\end{equation*}
for distinct \( i, j \in \{1, 2, \ldots, n\} \) such that \( \{i,j\} \in E \), is exactly the 0-th sheaf Laplacian matrix \( \mathbf{L}_{\mathcal{F}}^0 \).
\end{exam.}
\begin{proof}
To verify that the ANM Hessian matrix is exactly the sheaf Laplacian \( \mathbf{L}_{\mathcal{F}}^0 \) of the sheaf defined in this example, see~\cite{hu2024physicsinformedsheafmodel} for a detailed proof.   
\end{proof}

\paragraph{Tensor product of cellular sheaves}
At the end of Section~\ref{Section: Cellular Sheaves on Simplicial Complexes}, we introduce a sheaf operation called the \textit{tensor product}, which will be frequently used to construct new sheaf structures from a given collection of cellular sheaves on a simplicial complex. Specifically, let \( (K, \leq) \) be a simplicial complex, and let $\mathcal{F}_1, ..., \mathcal{F}_N: (K, \leq) \rightarrow \textup{\textsf{Mod}}_R$ be cellular sheaves over $K$. The \textit{tensor product} of these sheaves, denoted as $\mathcal{F}_1 \otimes \cdots \otimes \mathcal{F}_N$, is a cellular sheaf over $K$ defined as follows:
\begin{itemize}
\item[\rm (a)] $(\mathcal{F}_1 \otimes \cdots \otimes \mathcal{F}_N)_\sigma = (\mathcal{F}_1)_\sigma \otimes \cdots \otimes (\mathcal{F}_N)_\sigma$ whenever $\sigma \in K$;
\item[\rm (b)] $(\mathcal{F}_1 \otimes \cdots \otimes \mathcal{F}_N)_{\sigma,\tau} = (\mathcal{F}_1)_{\sigma,\tau} \otimes \cdots \otimes (\mathcal{F}_N)_{\sigma,\tau}$ whenever $\sigma \leq \tau$ in $K$.
\end{itemize}
Because the tensor product of homomorphisms automatically preserves the commutativity of the restriction maps, the assignment \( \mathcal{F}_1 \otimes \cdots \otimes \mathcal{F}_N : (K, \leq) \rightarrow \textup{\textsf{Mod}}_R \) naturally defines a cellular sheaf.

In Section~\ref{Section: Main-result-2}, we demonstrate how the tensor product of sheaves can be used to construct new sheaf structures on an underlying simplicial complex with higher-dimensional simplices. From this tensor product perspective, the sheaf structure motivated by the ANM model allows the vertex-edge, edge-face, and higher-dimensional relationships to be defined separately and then combined, providing a systematic approach to constructing sheaves on higher-dimensional objects. Moreover, as shown in Example~\ref{Example: main result-2}, the extended sheaf retains the Hessian information while additionally encoding higher-order interactions between simplices through higher-dimensional sheaf Laplacians.

\section{Cellular Sheaves on Higher-Dimensional Simplicial Complexes}
\label{Section: Cellular Sheaves on Higher-Dimensional Simplicial Complexes}

This section develops methods for constructing non-trivial cellular sheaves on simplicial complexes of dimension $\geq 2$. Sections \ref{Section: Main-result-1} and \ref{Section: Main-result-2} extend the ANM-based framework by incorporating higher-order geometric interactions through direct generalization and tensor product constructions, preserving the physical intuition from molecular modeling. In contrast, Section \ref{Section: Main-result-2} introduces an algebraically grounded approach using tools from commutative algebra, such as ringed spaces and sheaves of ideals or modules, without relying on any geometric or physical background. These two paradigms, the geometric-physical and the algebraic, offer complementary perspectives for constructing and interpreting cellular sheaf models on complex structures.

\subsection{Direct Extension of the ANM-Based Sheaf}\label{Section: Main-result-1}

The first method proposed in this paper involves a direct construction and extension of the ANM-based sheaf structure introduced in Example~\ref{Example: ANM sheaf on Graphs}. In this approach, we focus primarily on the geometric relationships between edges and faces within a simplicial complex \( (K, \leq) \) of dimension 2. A non-trivial sheaf \( \mathcal{F} : (K, \leq) \rightarrow \textup{\textsf{vect}}_{\mathbb{R}} \) is constructed such that \( \mathcal{F}_{e,f} \circ \mathcal{F}_{v,e} = 0 \) whenever \( v \in K_{(0)} \), \( e \in K_{(1)} \), and \( f \in K_{(2)} \). In particular, this construction satisfies the compositional axiom for restriction maps.

\begin{exam.}\label{Example: Direct definition-1}
Let \( (K, \leq) \) be a simplicial complex of dimension \( 2 \) with vertex set \( V = \{ 1, 2, \ldots, n \} \), where each vertex \( i \) is assigned a column vector \( \mathbf{r}_i \in \mathbb{R}^3 \). Furthermore, every edge $[i,j] \in K_{(1)}$ is associated with a weight $w_{ij} \in \mathbb{R}$, and every \( 2 \)-simplex \( \sigma = [i, j, k] \) is associated with a vector \( \mathbf{v}_{ijk} \) that is perpendicular to the affine space spanned by the vectors $ \mathbf{r}_i, \mathbf{r}_j, \mathbf{r}_k$. We define a cellular sheaf \( \mathcal{F} : (K, \leq) \rightarrow \textup{\textsf{vect}}_{\mathbb{R}} \) as follows:
\begin{itemize}
\item[\rm (a)] $\mathcal{F}_i = \mathbb{R}$ for every vertex $i \in V$.
\item[\rm (b)] $\mathcal{F}_{[i,j]} = \mathbb{R}^3$ for every $1$-simplex $[i,j] \in K_{(1)}$.
\item[\rm (c)] $\mathcal{F}_{[i,j,k]} = \mathbb{R}$ for every $2$-simplex $[i,j,k] \in K_{(2)}$.
\item[\rm (d)] \( \mathcal{F}_{i,[i,j]} = \mathcal{F}_{j,[i,j]} = w_{ij} \cdot (\mathbf{r}_j - \mathbf{r}_i) \), a \( 3 \times 1 \) matrix interpreted as an \( \mathbb{R} \)-linear map from \( \mathbb{R} \) to \( \mathbb{R}^3 \).
\item[\rm (e)] \( \mathcal{F}_{[i,j],[i,j,k]} = \mathcal{F}_{[i,k],[i,j,k]} = \mathcal{F}_{[j,k],[i,j,k]} = \mathbf{v}_{ijk}^\top \), a \( 1 \times 3 \) matrix interpreted as an \( \mathbb{R} \)-linear map from \( \mathbb{R}^3 \) to \( \mathbb{R} \).
\item[\rm (f)] $\mathcal{F}_{i,[i,j,k]} = \mathcal{F}_{j,[i,j,k]} = \mathcal{F}_{k,[i,j,k]}$ is defined as the zero map for every $2$-simplex $[i,j,k] \in K_{(2)}$.
\end{itemize}
Then, the rules~\textup{(a)--(f)} define a cellular sheaf \( \mathcal{F} : (K, \leq) \rightarrow \textup{\textsf{vect}}_{\mathbb{R}} \) of \( \mathbb{F} \)-vector spaces over the simplicial complex \( K \). 
\end{exam.} 
\begin{proof}
We first note that the vector \( \mathbf{v}_{ijk} \) always exists, as it can be chosen as the cross product of any two of the vectors \( \mathbf{r}_k - \mathbf{r}_i \), \( \mathbf{r}_k - \mathbf{r}_j \), and \( \mathbf{r}_j - \mathbf{r}_i \). To show that $\mathcal{F}$ satisfies the conditions of sheaves, it is sufficient to check that the restriction maps satisfy the commutative diagram:
\begin{equation*}
\xymatrix@+1.0em{
                & \mathcal{F}_\sigma
                \ar[r]^{\mathcal{F}_{\sigma,\tau_1}}
                \ar[d]_{\mathcal{F}_{\sigma,\tau_2}}
				& \mathcal{F}_{\tau_1}
                \ar[d]^{\mathcal{F}_{\tau_1,\eta}}
                \\
                & \mathcal{F}_{\tau_2}
                \ar[r]_{\mathcal{F}_{\tau_2,\eta}}
                & \mathcal{F}_\eta
                \\
}
\end{equation*}
whenever $\sigma \in V = K_{(0)}$, $\tau_1, \tau_2 \in K_{(1)}$, and $\eta \in K_{(2)}$.  Suppose $\sigma = i$ and $\eta = [i,j,k]$, then there are exactly two routes from $\sigma$ to $\eta$; namely, $v_i \leq \tau_1 = [i,j] \leq [i,j,k]$ and $i \leq \tau_2 = [i,k] \leq [i,j,k]$. For every $r \in \mathbb{R} = \mathcal{F}_{i}$, we have
\begin{equation*}
\begin{split}
(\mathcal{F}_{\tau_1, \eta} \circ \mathcal{F}_{\sigma, \tau_1})(r) &= (rw_{ij}) \cdot  \mathbf{v}_{ijk}^\top \cdot (\mathbf{r}_j - \mathbf{r}_i) = 0 = (rw_{ik}) \cdot  \mathbf{v}_{ijk}^\top \cdot (\mathbf{r}_k - \mathbf{r}_i) = (\mathcal{F}_{\tau_2, \eta} \circ \mathcal{F}_{\sigma, \tau_2})(r).
\end{split}
\end{equation*}
This shows that the composition map $\mathcal{F}_{\sigma,\eta}$ is always zero for $\sigma \in K_{(0)}$ and $\eta \in K_{(2)}$; in particular, the rule $\mathcal{F}$ satisfies all the sheaf conditions and hence forms a cellular sheaf of $\mathbb{F}$-vector spaces over $(K,\leq)$. 
\end{proof}
Furthermore, we sculpt the cochain complex of $\mathcal{F}$ as follows. Let $C^q(K;\mathcal{F})$ be the $q$-th cochain space of the sheaf $\mathcal{F}$ with $q = 0, 1, 2$. Suppose $K_{(1)} = \{ e_1, e_2, ..., e_m \}$ and $K_{(2)} = \{ f_1, f_2, ..., f_l \}$, then the $0$-th coboundary matrix $\mathbf{C}^0: C^0(K;\mathcal{F}) \rightarrow C^1(K;\mathcal{F})$ is defined as the $3m \times n$ block matrix
\begin{equation*}
\mathbf{C}^0 = \begin{bmatrix}
[1:e_1] \cdot \mathcal{F}_{1, e_1} & [2:e_1] \cdot \mathcal{F}_{2, e_1} & \cdots & [n:e_1] \cdot \mathcal{F}_{n, e_1}\\
[1:e_2] \cdot \mathcal{F}_{1, e_2} & [2:e_2] \cdot \mathcal{F}_{2, e_2} & \cdots & [n:e_2] \cdot \mathcal{F}_{n, e_2}\\
\vdots & \vdots & \ddots & \vdots \\
[1:e_m] \cdot \mathcal{F}_{1, e_m} & [2:e_m] \cdot \mathcal{F}_{2, e_m} & \cdots & [n:e_m] \cdot \mathcal{F}_{n, e_m}\\
\end{bmatrix},    
\end{equation*}
where each row of blocks in \( \mathbf{C}^0 \) contains exactly two non-zero blocks, corresponding to the two endpoints of the edge associated with that row. On the other hand, the $1$-st coboundary matrix $\mathbf{C}^1: C^1(K;\mathcal{F}) \rightarrow C^2(K;\mathcal{F})$ is defined as the $l \times 3m$ block matrix
\begin{equation*}
\mathbf{C}^1 = \begin{bmatrix}
[e_1:f_1] \cdot \mathcal{F}_{e_1, f_1} & [e_2:f_1] \cdot \mathcal{F}_{e_2, f_1} & \cdots & [e_m:f_1] \cdot \mathcal{F}_{e_m, f_1}\\
[e_1:f_2] \cdot \mathcal{F}_{e_1, f_2} & [e_2:f_2] \cdot \mathcal{F}_{e_2, f_2} & \cdots & [e_m:f_2] \cdot \mathcal{F}_{e_m, f_2}\\
\vdots & \vdots & \ddots & \vdots \\
[e_1:f_m] \cdot \mathcal{F}_{e_1, f_l} & [e_2:f_l] \cdot \mathcal{F}_{e_2, f_l} & \cdots & [e_m:f_l] \cdot \mathcal{F}_{e_m, f_l}\\
\end{bmatrix},    
\end{equation*}
where each row of blocks in $\mathbf{C}^1$ only contains three non-zero blocks, as the corresponding $2$-simplex only involves $3$ $1$-faces. Therefore, the $0$-th, $1$-th, and $2$-nd sheaf Laplacians are therefore defined as
\begin{equation*}
\mathbf{L}_{\mathcal{F}}^0 = (\mathbf{C}^0)^{\top} \cdot \mathbf{C}^0, \ \mathbf{L}_{\mathcal{F}}^1 = \mathbf{C}^0 \cdot (\mathbf{C}^0)^{\top} + (\mathbf{C}^1)^{\top} \cdot \mathbf{C}^1, \text{ and } \mathbf{L}_{\mathcal{F}}^2 = \mathbf{C}^1 \cdot (\mathbf{C}^1)^{\top}.
\end{equation*}
In particular, by utilizing the block multiplication of block matrices, as a part of, the $n \times n$ $0$-th Laplacian matrix of $\mathcal{F}$ can be represented by the following form
\begin{equation*}
\mathbf{L}_{\mathcal{F}}^0 = \begin{bmatrix}
L_{11}^0 & L_{12}^0 & \cdots & L_{1n}^0\\
L_{21}^0 & L_{22}^0 & \cdots & L_{2n}^0\\
\vdots & \vdots & \ddots & \vdots \\
L_{n1}^0 & L_{n2}^0 & \cdots & L_{nn}^0\\
\end{bmatrix}, 
\end{equation*}
where $L_{ij}^0 = \sum_{r = 1}^m [i : e_r] \cdot [j : e_r] \cdot \mathcal{F}_{i, e_r}^\top \cdot \mathcal{F}_{j, e_r} = -w_{ij}^2 \cdot \Vert \mathbf{r}_j - \mathbf{r}_i \Vert^2$ for distinct $i, j \in \{ 1, 2, ..., n \}$ such that \( [i,j] \in K_{(1)} \), and $L_{ii}^0 = \sum_{j \neq i} - L_{ij}^0$ for $i \in \{ 1, 2, ..., n \}$. In other words, the $0$-th sheaf Laplacian of $\mathcal{F}$ collects all the distance information of edges within the underlying graph. 

On the other hand, as a part of the $3m \times 3m$ $1$-th Laplacian matrix $\mathbf{L}_{\mathcal{F}}^1$, the matrix $\mathbf{C}^0 \cdot (\mathbf{C}^0)^{\top}$ can be represented by the following block form
\begin{equation*}
\mathbf{C}^0 \cdot (\mathbf{C}^0)^{\top} = \begin{bmatrix}
L_{11} & L_{12} & \cdots & L_{1m}\\
L_{21} & L_{22} & \cdots & L_{2m}\\
\vdots & \vdots & \ddots & \vdots \\
L_{m1} & L_{m2} & \cdots & L_{mm}\\
\end{bmatrix}, 
\end{equation*}
where $L_{ab} = \sum_{k = 1}^n [k: e_a] \cdot [k: e_b] \cdot \mathcal{F}_{k,e_a} \cdot \mathcal{F}_{k,e_b}^{\top}$ for arbitrary $a, b \in \{ 1, 2, ..., m \}$. On the other hand, for every \( a \in \{ 1, 2, \ldots, m \} \), where \( e_a = [i,j] \) consists of two vertices \( i \) and \( j \), the block matrix \( L_{aa} \) is given by
\begin{equation*}
\begin{split}
2 \cdot \mathcal{F}_{i,[i,j]} \cdot \mathcal{F}_{i,[i,j]}^{\top} &= 2 \cdot w_{ij}^2 \cdot (\mathbf{r}_j - \mathbf{r}_i) \cdot (\mathbf{r}_j - \mathbf{r}_i)^{\top}
\\ 
&= 2 \cdot w_{ij}^2 \cdot \begin{bmatrix}
x_j^\circ - x_i^\circ \\
y_j^\circ - y_i^\circ \\
z_j^\circ - z_i^\circ \\
\end{bmatrix} \cdot \begin{bmatrix}
x_j^\circ - x_i^\circ & y_j^\circ - y_i^\circ & z_j^\circ - z_i^\circ \\
\end{bmatrix} \\
&= 2 \cdot w_{ij}^2 \cdot \begin{bmatrix}
(x_j^\circ - x_i^\circ)^2 & (x_j^\circ - x_i^\circ)(y_j^\circ - y_i^\circ) & (x_j^\circ - x_i^\circ)(z_j^\circ - z_i^\circ)\\
(x_j^\circ - x_i^\circ)(y_j^\circ - y_i^\circ) & (y_j^\circ - y_i^\circ)^2 & (y_j^\circ - y_i^\circ)(z_j^\circ - z_i^\circ)\\
(x_j^\circ - x_i^\circ)(z_j^\circ - z_i^\circ) & (y_j^\circ - y_i^\circ)(z_j^\circ - z_i^\circ) & (z_j^\circ - z_i^\circ)^2\\
\end{bmatrix},
\end{split}
\end{equation*}
where $\mathbf{r}_i = (x_i^\circ, y_i^\circ, z_i^\circ)$. In other words, the original \( 3n \times 3n \) Hessian information is, in a certain sense, encoded within the diagonal blocks of the matrix \( L \), since each off-diagonal block \( H_{ij} \) (with \( i \neq j \)) in the ANM Hessian matrix \( \mathbf{H}_{\operatorname{ANM}} \) appears as a component within the diagonal blocks of \( \mathbf{C}^0 \cdot (\mathbf{C}^0)^{\top} \). On the other hand, for $a \neq b$, say $e_a = [i, j]$ and $e_b = [i, k]$ (or, $e_a = [j, i]$ and $e_b = [k, i]$). Then, 
\begin{equation*}
\begin{split}
L_{ab} &= (w_{ij} w_{ik}) \cdot \mathcal{F}_{i,[i,j]} \cdot \mathcal{F}_{i,[i,k]}^{\top} \\
&= (w_{ij} w_{ik}) \cdot \begin{bmatrix}
(x_j^\circ - x_i^\circ)(x_k^\circ - x_i^\circ) & (x_j^\circ - x_i^\circ)(y_k^\circ - y_i^\circ) & (x_j^\circ - x_i^\circ)(z_k^\circ - z_i^\circ)\\
(y_j^\circ - y_i^\circ)(x_k^\circ - x_i^\circ) & (y_j^\circ - y_i^\circ)(y_k^\circ - y_i^\circ) & (y_j^\circ - y_i^\circ)(z_k^\circ - z_i^\circ)\\
(z_j^\circ - z_i^\circ)(x_k^\circ - x_i^\circ) & (z_j^\circ - z_i^\circ)(y_k^\circ - y_i^\circ) & (z_j^\circ - z_i^\circ)(z_k^\circ - z_i^\circ)\\
\end{bmatrix}.
\end{split}   
\end{equation*}
On the other hand, if $e_a = [i, j]$ and $e_b = [k, i]$, then, $L_{ab} = (-w_{ij} w_{ik}) \cdot \mathcal{F}_{i,[i,j]} \cdot \mathcal{F}_{i,[i,k]}^{\top}$. In particular, each submatrix within the off-diagonal blocks of \( \mathbf{C}^0 (\mathbf{C}^0)^{\top} \) encodes quadratic information about the coordinates of three interacting vertices connected by linkage relations.

Furthermore, the matrices \( (\mathbf{C}^1)^{\top} \mathbf{C}^1 \) and \( \mathbf{L}_{\mathcal{F}}^2 = \mathbf{C}^1 (\mathbf{C}^1)^{\top} \) encode more complex relationships between edges and faces, providing richer information for analyzing the structure of the simplicial complex.

\subsection{Tensor Product Extension of the ANM-Based Sheaf}\label{Section: Main-result-2}
The second method employs the tensor product of sheaves to construct new sheaf structures by combining finitely many sheaves, each capturing distinct relationships between simplices of different dimensions. As described in Section~\ref{Section: Main-result-1}, we focus on ANM-based sheaves defined on graphs embedded in \( \mathbb{R}^2 \), demonstrating that the original ANM Hessian information is retained within the extended sheaf defined on a simplicial complex containing higher-dimensional simplices. Further exploration of the geometric meaning of higher-dimensional sheaf Laplacians remains a promising direction for future research.

\begin{exam.}\label{Example: main result-2}
Let \( (K, \leq) \) be a simplicial complex of dimension \( 2 \) with vertex set \( V = \{ 1, 2, \ldots, n \} \), where each vertex \( i \) is assigned a column vector \( \mathbf{r}_i \in \mathbb{R}^3 \). Furthermore, each edge \( [i,j] \in K_{(1)} \) is associated with a weight \( w_{ij} \in \mathbb{R} \). Consider two cellular sheaves defined on \( K \), \( \mathcal{F}, \mathcal{G} : (K, \leq) \rightarrow \textup{\textsf{vect}}_{\mathbb{R}} \). Specifically, we define \( \mathcal{F} \) as the cellular sheaf in Example \ref{Example: Direct definition-1}.

On the other hand, the cellular sheaf $\mathcal{G}: (K, \leq) \rightarrow \textup{\textsf{vect}}_{\mathbb{R}}$ is defined as follows. For every \( \sigma \in K_{(2)} \), we associate a column vector \( \mathbf{w}_\sigma \in \mathbb{R}^3 \). We define a cellular sheaf \( \mathcal{G} : (K, \leq) \rightarrow \textup{\textsf{vect}}_{\mathbb{R}} \) by the following assignments:
\begin{itemize}
\item[\rm (a)] $\mathcal{G}_i = \mathbb{R}^3$ for every $i \in V$.
\item[\rm (b)] $\mathcal{G}_{[i,j]} = \mathbb{R}^3$ for every $[i,j] \in K_{(1)}$.
\item[\rm (c)] $\mathcal{G}_{[i,j,k]} = \mathbb{R}$ for every $[i,j,k] \in K_{(2)}$.
\item[\rm (d)] $\mathcal{G}_{i,[i,j]} = \mathcal{G}_{j,[i,j]}$ is defined as the identity map ${\rm id}_{\mathbb{R}^3}$.
\item[\rm (e)] For \( \sigma \in K_{(1)} \) and \( \tau \in K_{(2)} \) with \( \sigma \leq \tau \), the map \( \mathcal{G}_{\sigma, \tau} \) is defined as the \( 1 \times 3 \) matrix \( \mathbf{w}_\tau^\top \).
\end{itemize}
Because the maps $\mathcal{G}_{\bullet,\bullet}$ defined in (d) are identity maps, the rule $\mathcal{G}$ automatically forms a cellular sheaf over the simplicial complex $K$. Then, the tensor product $\mathcal{F} \otimes \mathcal{G}$ of cellular sheaves $\mathcal{F}$ and $\mathcal{G}$ is defined.
\end{exam.}
As presented in Section~\ref{Section: Main-result-1}, we show that the sheaf \( \mathcal{F} \otimes \mathcal{G} \) retains the information of the original ANM Hessian matrix associated with the cellular sheaf defined in Example~\ref{Example: ANM sheaf on Graphs}. The details of the sheaf \( \mathcal{F} \otimes \mathcal{G} \) are provided below.
\begin{itemize}
\item[\rm (a)] $(\mathcal{F} \otimes \mathcal{G})_i = \mathbb{R}^3 \otimes \mathbb{R} \simeq \mathbb{R}^3$ for every $i \in V$.
\item[\rm (b)] $(\mathcal{F} \otimes \mathcal{G})_{[i,j]} = \mathbb{R}^3 \otimes \mathbb{R}^3 \simeq \mathbb{R}^9$ for every $[i,j] \in K_{(1)}$.
\item[\rm (c)] $(\mathcal{F} \otimes \mathcal{G})_{[i,j,k]} = \mathbb{R} \otimes \mathbb{R} \simeq \mathbb{R}$ for every $[i,j,k] \in K_{(2)}$.
\end{itemize}
Suppose $K_{(1)} = \{ e_1, e_2, ..., e_m \}$, the $0$-coboundary matrix $\mathbf{C}_{\mathcal{F} \otimes \mathcal{G}}^0: C^0(K;\mathcal{F}\otimes\mathcal{G}) \rightarrow C^1(K;\mathcal{F}\otimes\mathcal{G})$ is
\begin{equation*}
\begin{bmatrix}
[1:e_1] \cdot (\mathcal{F} \otimes \mathcal{G})_{1, e_1} & \cdots & [n:e_1] \cdot (\mathcal{F} \otimes \mathcal{G})_{n, e_1}\\
[1:e_2] \cdot (\mathcal{F} \otimes \mathcal{G})_{1, e_2} & \cdots & [n:e_2] \cdot (\mathcal{F} \otimes \mathcal{G})_{n, e_2}\\
\vdots & \ddots & \vdots \\
[1:e_m] \cdot (\mathcal{F} \otimes \mathcal{G})_{1, e_m} & \cdots & [n:e_m] \cdot (\mathcal{F} \otimes \mathcal{G})_{n, e_m}\\
\end{bmatrix}.
\end{equation*}
Because $(\mathcal{F} \otimes \mathcal{G})_{k, e_r} = \mathcal{F}_{k, e_r} \otimes \mathcal{G}_{k, e_r} = \mathcal{F}_{k, e_r} \otimes I_3$ for every $k \in \{ 1, 2, ..., n \}$ and $r \in \{ 1, 2, ..., m \}$, the $0$-th cobounary matrix $\mathbf{C}_{\mathcal{F} \otimes \mathcal{G}}^0$ is exactly the $9m \times 3n$ matrix
\begin{equation}\label{Eq. 0-th sheaf coboundary map}
\mathbf{C}_{\mathcal{F} \otimes \mathcal{G}}^0 = \begin{bmatrix}
[1:e_1] \cdot \mathcal{F}_{1, e_1} & \cdots & [n:e_1] \cdot \mathcal{F}_{n, e_1}\\
[1:e_2] \cdot \mathcal{F}_{1, e_2} & \cdots & [n:e_2] \cdot \mathcal{F}_{n, e_2}\\
\vdots & \ddots & \vdots \\
[1:e_m] \cdot \mathcal{F}_{1, e_m} & \cdots & [n:e_m] \cdot \mathcal{F}_{n, e_m}\\
\end{bmatrix} \otimes \begin{bmatrix}
1 & 0 & 0 \\
0 & 1 & 0 \\
0 & 0 & 1 \\
\end{bmatrix} = \mathbf{C}_{\mathcal{F}}^0 \otimes \begin{bmatrix}
1 & 0 & 0 \\
0 & 1 & 0 \\
0 & 0 & 1 \\
\end{bmatrix},  
\end{equation}
where $\mathbf{C}_{\mathcal{F}}^0:  C^0(K;\mathcal{F}) \rightarrow C^1(K;\mathcal{F})$ is the $0$-th coboundary matrix of $\mathcal{F}$. Therefore, 
the $0$-th sheaf Laplacian of $\mathcal{F} \otimes \mathcal{G}$ can be expressed as the $3n \times 3n$ matrix:
\begin{equation}\label{Eq. 0-th sheaf Laplacian of the tensor product}
\begin{split}
\mathbf{L}_{\mathcal{F} \otimes \mathcal{G}}^0 &= (\mathbf{C}_{\mathcal{F} \otimes \mathcal{G}}^0)^{\top} \cdot (\mathbf{C}_{\mathcal{F} \otimes \mathcal{G}}^0) \\ 
&= \left( \mathbf{C}_{\mathcal{F}}^0 \otimes \begin{bmatrix}
1 & 0 & 0 \\
0 & 1 & 0 \\
0 & 0 & 1 \\
\end{bmatrix} \right)^{\top} \cdot \left( \mathbf{C}_{\mathcal{F}}^0 \otimes \begin{bmatrix}
1 & 0 & 0 \\
0 & 1 & 0 \\
0 & 0 & 1 \\
\end{bmatrix} \right) \\
&= \left( (\mathbf{C}_{\mathcal{F}}^0)^{\top} \otimes \begin{bmatrix}
1 & 0 & 0 \\
0 & 1 & 0 \\
0 & 0 & 1 \\
\end{bmatrix} \right) \cdot \left( \mathbf{C}_{\mathcal{F}}^0 \otimes \begin{bmatrix}
1 & 0 & 0 \\
0 & 1 & 0 \\
0 & 0 & 1 \\
\end{bmatrix} \right) \\
&= ((\mathbf{C}_{\mathcal{F}}^0)^{\top}\mathbf{C}_{\mathcal{F}}^0) \otimes \begin{bmatrix}
1 & 0 & 0 \\
0 & 1 & 0 \\
0 & 0 & 1 \\
\end{bmatrix} 
= \mathbf{L}_{\mathcal{F}}^0 \otimes \begin{bmatrix}
1 & 0 & 0 \\
0 & 1 & 0 \\
0 & 0 & 1 \\
\end{bmatrix}.
\end{split}     
\end{equation}
On the other hand,
\begin{equation}\label{Eq. Lower part of the 1-th sheaf Laplacian of the tensor product}
(\mathbf{C}_{\mathcal{F} \otimes \mathcal{G}}^0) \cdot (\mathbf{C}_{\mathcal{F} \otimes \mathcal{G}}^0)^{\top} =(\mathbf{C}_{\mathcal{F}}^0(\mathbf{C}_{\mathcal{F}}^0)^{\top}) \otimes \begin{bmatrix}
1 & 0 & 0 \\
0 & 1 & 0 \\
0 & 0 & 1 \\
\end{bmatrix}.
\end{equation}
In other words, by the result of Example \ref{Example: Direct definition-1}, the lower part $\mathbf{C}_{\mathcal{F}}^0(\mathbf{C}_{\mathcal{F}}^0)$ of the $1$-th sheaf Laplacian of \( \mathcal{F} \) preserves the original Hessian information in its diagonal blocks, and \( (\mathbf{C}_{\mathcal{F} \otimes \mathcal{G}}^0) \) extends it across three dimensions. Alternatively, by defining meaningful linear transformations \( \mathcal{G}_{\sigma, \tau} \) for \( \sigma \in K_{(1)} \) and \( \tau \in K_{(2)} \) satisfying \( \sigma \leq \tau \), the higher-dimensional sheaf Laplacian of \( \mathcal{F} \otimes \mathcal{G} \) can encode the edge-face relationships within the simplicial complex. This construction entangles relationships across different dimensions within \( \mathbf{L}^1_{\mathcal{F} \otimes \mathcal{G}} = \mathbf{C}_{\mathcal{F} \otimes \mathcal{G}}^0 \cdot (\mathbf{C}_{\mathcal{F} \otimes \mathcal{G}}^0)^{\top} + (\mathbf{C}_{\mathcal{F} \otimes \mathcal{G}}^1)^{\top} \cdot \mathbf{C}_{\mathcal{F} \otimes \mathcal{G}}^1\).

\begin{remark}
In the construction above, the restriction maps \( \mathcal{G}_{i,[i,j]} = \mathcal{G}_{j,[i,j]} \) are defined as identity maps. However, various alternative choices can also be considered. For instance, if all restriction maps \( \mathcal{G}_{i,\sigma} \) for \( \sigma \in K_{(2)} \) are assigned the same linear transformations, the composition rules for restriction maps still hold, and the resulting assignment remains a valid cellular sheaf. Under this modified setting, the corresponding coboundary maps and sheaf Laplacians, such as those defined in \eqref{Eq. 0-th sheaf coboundary map}, \eqref{Eq. 0-th sheaf Laplacian of the tensor product}, and \eqref{Eq. Lower part of the 1-th sheaf Laplacian of the tensor product}, will change accordingly. These changes may lead to distinct algebraic structures and geometric interpretations.
\end{remark}

\subsection{Cellular Sheaves of Rings and Modules}\label{Section: Main-result-3}

The third approach explores general cellular sheaf representations on simplicial complexes with higher-dimensional simplices by employing tools from commutative algebra, without relying on a specific geometric or physical background as in Sections~\ref{Section: Main-result-1} and~\ref{Section: Main-result-2}. In particular, by considering the structure of a ringed space and sheaves of \( R \)-modules over a given commutative ring \( R \), natural sheaf structures can be defined in a straightforward manner.

\paragraph{Ringed Spaces on Simplicial Complexes}

Ringed spaces on a topological space play a central role in algebraic and differential geometry, appearing in fundamental constructions such as schemes and sheaves of \( \mathscr{C}^\infty \) functions on differentiable manifolds~\cite{hartshorne2013algebraic,bredon2012sheaf}. For a comprehensive introduction to the commutative algebra tools used in this paper, see~\cite{Atiyah, matsumura1989commutative}.

\begin{def.}
A \textbf{ringed space} is a pair \( (X, \mathcal{O}_X) \), where \( X \) is a topological space and \( \mathcal{O}_X \) is a sheaf of rings on \( X \).
\end{def.}

In discrete settings, such as cellular sheaves on a simplicial complex \( K \) equipped with its face relations, a \textit{ringed space} with base space \( K \) is given by the pair \( (K, \mathcal{O}) \),  where \( \mathcal{O} : (K, \leq) \rightarrow \textup{\textsf{Ring}} \) is a cellular sheaf of rings on \( K \). In fact, when \( K \) is endowed with the Alexandrov topology \( \mathfrak{A} \) induced by the partial order \( \leq \)~\cite{Ale37,Ale47}, any functor from \( (K, \leq) \) to \( \textup{\textsf{Ring}} \) (or to \( \textup{\textsf{Mod}}_R, \textup{\textsf{Vect}}_{\mathbb{F}} \), or \( \textsf{Ab} \)) defines a sheaf on the topological space \( (K, \mathfrak{A}) \); see, for example,~\cite{curry2014sheaves,hu2020briefnotesheafstructures,ayzenberg2025sheaf}.

As one of the most important algebraic structures associated with a cellular sheaf, determining the ring of global sections of a ringed space is a fundamental task. The following proposition provides a straightforward yet essential identification.

\begin{prop.}\label{Proposition: the ring of global sections}
Let \( (K, \mathcal{O}) \) be a ringed space based on a simplicial complex $K$. Then, the ring of the global sections of $\mathcal{O}$ is the inverse limit
\begin{equation*}
\Gamma(K;\mathcal{O}) = \varprojlim_{\sigma \in (K, \leq)} \mathcal{O}_\sigma = \left\{ \left. (r_\sigma)_{\sigma \in K} \in \prod_{\sigma \in K} \mathcal{O}_\sigma \ \right| \ \mathcal{O}_{\sigma, \eta}(r_\sigma) = \mathcal{O}_{\tau, \eta}(r_\tau) \text{ if } \sigma \leq \eta \text{ and } \tau \leq \eta \right\}.    
\end{equation*}
\end{prop.}
\begin{proof}
This proposition is an immediate consequence of Equation~\eqref{Eq. global section space as an inverse limit}. Note that \( \Gamma(K; \mathcal{O}) \) contains a multiplicative identity element, since each ring \( \mathcal{O}_\sigma \) contains the identity element \( 1_{\mathcal{O}_\sigma} \), and every ring homomorphism \( \mathcal{O}_{\sigma, \tau} \) with \( \sigma \leq \tau \) preserves identity elements. Therefore, the constant section \( (1_{\mathcal{O}_\sigma})_{\sigma \in K} \) belongs to \( \Gamma(K; \mathcal{O}) \).
\end{proof}

In particular, as illustrated in the following example, for rings \( A, B, C \) and ring homomorphisms \( f : A \rightarrow C \) and \( g : B \rightarrow C \), the \textit{fibre product} of \( A \) and \( B \) over \( C \) through the homomorphisms $f$ and $g$, denoted by $A \times_C B := \{ (a, b) \in A \times B \mid f(a) = g(b) \}$, can be viewed as a special case of the ring of global sections of a cellular sheaf over the simplicial complex \( K = \{ \{1\}, \{2\}, \{1, 2\} \} \).

\begin{exam.}[Fibre product of rings]
Let \( K \) be the simplicial complex \( \{ \{1\}, \{2\}, \{1, 2\} \} \). Consider the cellular sheaf \( \mathcal{O} : (K, \leq) \rightarrow \textup{\textsf{Ring}} \) defined as follows:
\begin{itemize}
\item[\rm (a)] $\mathcal{O}_{\{1\}} = A$, $\mathcal{O}_{\{2\}} = B$, and $\mathcal{O}_{\{1,2\}} = C$;
\item[\rm (a)] $\mathcal{O}_{\{1\},\{1,2\}} = f$, $\mathcal{O}_{\{2\},\{1,2\}} = g$, and $\mathcal{O}_{\sigma\sigma} = \id_{\mathcal{O}_\sigma}$ for $\sigma \in K$.
\end{itemize}
Then the ring of global sections $\Gamma(K, \mathcal{O})$ is isomorphic to $A \times_C B$ as rings.
\end{exam.} 
\begin{proof}
The cellular sheaf \( \mathcal{O} : (K, \leq) \rightarrow \textup{\textsf{Ring}} \) can be represented by the following diagram:
\begin{equation*}
\xymatrix@+2.0em{
                & 
				& \mathcal{O}_{\{ 1 \}}
                \ar[d]^{\mathcal{O}_{\{ 1 \},\{ 1, 2 \}}}
                \\
                & \mathcal{O}_{\{ 2 \}}
                \ar[r]_{\mathcal{O}_{\{ 2 \},\{ 1, 2 \}}}
                & \mathcal{O}_{\{ 1, 2 \}}
                \\
}
\end{equation*}
By \eqref{Eq. global section space as an inverse limit-simple form}, $\Gamma(K, \mathcal{O}) \simeq \{ (a, b) \in \mathcal{O}_{\{ 1 \}} \times \mathcal{O}_{\{ 2 \}} \mid \mathcal{O}_{\{ 1 \},\{ 1, 2 \}}(a) = \mathcal{O}_{\{ 2 \},\{ 1, 2 \}}(b) \}$. By the definition of the cellular sheaf $\mathcal{O}$, the ring in the right-hand side is exact the fibre product $\{\{ (a, b) \in A \times B \mid f(a) = g(b) \} = A \times_C B$.
\end{proof} 

Focusing on assigning rings to the simplices of a simplicial complex, we present several straightforward methods for constructing cellular sheaves on simplicial complexes.

\begin{prop.}\label{Proposition: sheaf of rings induced by a sheaf of ideals}
Let \( R \) be a ring and \( K \) a simplicial complex with face relation \( \leq \). Let \( I(R) \) denote the set of ideals of \( R \). Viewing \( (I(R), \subseteq) \) as a poset category, any functor \( \mathcal{I} : (K, \leq) \rightarrow (I(R), \subseteq) \) defines a cellular sheaf of rings \( \mathcal{O} : (K, \leq) \rightarrow \textup{\textsf{Ring}} \) by setting \( \mathcal{O}_\sigma = R / \mathcal{I}_\sigma \) for each \( \sigma \in K \), and defining the restriction map \( \mathcal{O}_{\sigma, \tau} : R / \mathcal{I}_\sigma \rightarrow R / \mathcal{I}_\tau \) for \( \sigma \leq \tau \) as the canonical ring homomorphism \( \overline{x} \mapsto \overline{x} \).
\end{prop.}
\begin{proof}
Because \( \mathcal{I} : (K, \leq) \rightarrow (I(R), \subseteq) \) is a functor, $\mathcal{I}_{\sigma} \subseteq \mathcal{I}_{\tau}$ whenever $\sigma \leq \tau$. This shows that the map \( \mathcal{O}_{\sigma, \tau} : R / \mathcal{I}_\sigma \rightarrow R / \mathcal{I}_\tau \) defined by \( \overline{x} \mapsto \overline{x} \) is a well-defined ring homomorphism. Clearly, $\mathcal{O}_{\sigma, \rho} = \mathcal{O}_{\tau, \rho} \circ \mathcal{O}_{\sigma, \tau}$ for $\sigma \leq \tau \leq \rho$ in $K$, and this shows that $\mathcal{O}$ forms a cellular sheaf of rings on the simplicial complex $K$.
\end{proof}  

Using the language of sheaf theory, the functor \( \mathcal{I} : (K, \leq) \rightarrow (I(R), \subseteq) \) can be interpreted as a \textit{cellular sheaf of ideals} in a ring \( R \). In particular, the sheaf \( \mathcal{O} \) constructed in Proposition~\ref{Proposition: sheaf of rings induced by a sheaf of ideals} is precisely the \textit{quotient sheaf} \( \underline{R} / \mathcal{I} \), where \( \underline{R} \) denotes the constant sheaf of rings with value \( R \) on \( K \).

In addition to the quotient construction in Proposition~\ref{Proposition: sheaf of rings induced by a sheaf of ideals}, another approach for constructing cellular sheaves on simplicial complexes involves the use of rings of fractions. Recall that for a ring $A$ and a prime ideal $\mathfrak{p} \in \operatorname{Spec}(A)$, the localization of $A$ at prime ideal $\mathfrak{p}$ is denoted by $A_{\mathfrak{p}} = \{ a/s \mid a \in A, s \in A \setminus S \}$.
 
\begin{prop.}\label{Proposition: sheaf of rings induced by localization}
Let \( R \) be a ring and \( K \) a simplicial complex with face relation \( \leq \). Let \( \operatorname{Spec}(R) \) denote the set of prime ideals of \( R \). Viewing \( (\operatorname{Spec}(R), \subseteq) \) as a poset category, any functor \( \mathcal{P} : (K, \leq)^{\operatorname{op}} \rightarrow (\operatorname{Spec}(R), \subseteq) \) defines a cellular sheaf of rings \( \mathcal{O} : (K, \leq) \rightarrow \textup{\textsf{Ring}} \) by setting \( \mathcal{O}_\sigma = R_{\mathcal{P}_\sigma} \) for each \( \sigma \in K \), and defining the restriction map \( \mathcal{O}_{\sigma, \tau} : R_{\mathcal{P}_\sigma} \rightarrow R_{\mathcal{P}_\tau} \) for \( \sigma \leq \tau \) as the canonical ring homomorphism \( r/s \mapsto r/s \).
\end{prop.}
\begin{proof}
Because \( \mathcal{P} : (K, \leq)^{\operatorname{op}} \rightarrow (\operatorname{Spec}(R), \subseteq) \) is a functor, $\mathcal{P}_{\tau} \subseteq \mathcal{P}_{\sigma}$ whenever $\sigma \leq \tau$. By the universal property of rings of fractions, there is a canonical ring homomorphism $\mathcal{O}_{\sigma, \tau}: R_{\mathcal{P}_\sigma} \rightarrow R_{\mathcal{P}_\tau}$ so that $\mathcal{O}_{\sigma, \tau}(r/s) = r/s$ for every $r/s \in R_{\mathcal{P}_\sigma}$ with $r \in R$ and $s \in R \setminus \mathcal{P}_\sigma$. Clearly, $\mathcal{O}_{\sigma, \rho} = \mathcal{O}_{\tau, \rho} \circ \mathcal{O}_{\sigma, \tau}$ for $\sigma \leq \tau \leq \rho$ in $K$, and this shows that $\mathcal{O}$ forms a cellular sheaf of rings on the simplicial complex $K$.
\end{proof}

Corresponding to Propositions~\ref{Proposition: sheaf of rings induced by a sheaf of ideals} and~\ref{Proposition: sheaf of rings induced by localization}, we present the following examples of cellular sheaves of rings on a simplicial complex, each based on constructions in polynomial rings.

\begin{exam.}\label{Example: Quotient sheaf-1}
Let \( (K, \leq) \) be a simplicial complex with vertex set \( V = \{ 1, 2, \ldots, n \} = K_{(0)} \). Let \( \mathbb{F} \) be a field, and let \( \mathbb{S} = \mathbb{F}[x_1, \ldots, x_n] \) denote the polynomial ring in \( n \) variables over \( \mathbb{F} \). For each simplex \( \sigma = \{ i_0, i_1, \ldots, i_q \} \subseteq \{ 1, 2, \ldots, n \} \), we define \( \mathcal{I}_\sigma \) to be the ideal \( (x_{i_0}, x_{i_1}, \ldots, x_{i_q}) \subseteq \mathbb{S} \). Then \( \mathcal{I}_\sigma \subseteq \mathcal{I}_\tau \) whenever \( \sigma \leq \tau \) in \( K \). This assignment defines a cellular sheaf \( \mathcal{I} : (K, \leq) \rightarrow (I(\mathbb{S}), \subseteq) \). In particular, the quotient cellular sheaf $\underline{\mathbb{S}}/\mathcal{I}$ forms a cellular sheaf of rings on the simplicial complex $K$.
\end{exam.}

The cellular sheaf constructed in Example~\ref{Example: Quotient sheaf-1} is based on the assignment where each ideal \( \mathcal{I}_\sigma \) is a prime ideal. The following example illustrates a construction based on a different type of ideals.

\begin{exam.}
Let \( (K, \leq) \) be a simplicial complex with vertex set \( V = \{ 1, 2, \ldots, n \} = K_{(0)} \). Let \( \mathbb{F} \) be a field, and let \( \mathbb{S} = \mathbb{F}[x_1, \ldots, x_n] \) denote the polynomial ring in \( n \) variables over \( \mathbb{F} \). For each simplex \( \sigma = \{ i_0, i_1, \ldots, i_q \} \subseteq \{ 1, 2, \ldots, n \} \), we define 
\begin{equation*}
\mathcal{I}_\sigma = (x_{i}x_j \mid \{ i,j \} \subseteq \sigma ).    
\end{equation*}
Then \( \mathcal{I}_\sigma \subseteq \mathcal{I}_\tau \) whenever \( \sigma \leq \tau \) in \( K \). This assignment defines a cellular sheaf \( \mathcal{I} : (K, \leq) \rightarrow (I(\mathbb{S}), \subseteq) \). In particular, the quotient cellular sheaf $\underline{\mathbb{S}}/\mathcal{I}$ forms a cellular sheaf of rings on the simplicial complex $K$.    
\end{exam.}

\begin{exam.}
Let \( (K, \leq) \) be a simplicial complex with vertex set \( V = \{ 1, 2, \ldots, n \} = K_{(0)} \). Let \( \mathbb{F} \) be a field, and let \( \mathbb{S} = \mathbb{F}[x_1, \ldots, x_n] \) denote the polynomial ring in \( n \) variables over \( \mathbb{F} \). For each simplex \( \sigma = \{ i_0, i_1, \ldots, i_q \} \subseteq \{ 1, 2, \ldots, n \} \), we define \( \mathcal{I}_\sigma \) to be the prime ideal 
\begin{equation*}
\mathcal{P}_\sigma = (x_{i} \mid i \notin \sigma) \subseteq \mathbb{S}.    
\end{equation*}
Then \( \mathcal{P}_\tau \subseteq \mathcal{P}_\sigma \) whenever \( \sigma \leq \tau \) in \( K \). This assignment defines a cellular sheaf \( \mathcal{P} : (K, \leq)^{\operatorname{op}} \rightarrow (\operatorname{Spec}(\mathbb{S}), \subseteq) \). In particular, by Proposition \ref{Proposition: sheaf of rings induced by localization}, it induces a cellular sheaf \( \mathcal{O} : (K, \leq) \rightarrow \textup{\textsf{Ring}} \) so that $\mathcal{O}_\sigma = R_{\mathcal{P}_\sigma}$ whenever $\sigma \in K$.
\end{exam.}

With different purposes and approaches, various cellular sheaves \( \mathcal{I} : (K, \leq) \rightarrow (I(\mathbb{S}), \subseteq) \) and \( \mathcal{P} : (K, \leq)^{\operatorname{op}} \rightarrow (\operatorname{Spec}(\mathbb{S}), \subseteq) \) can be defined based on the combinatorial properties of the underlying simplicial complex. Exploring the connections between these induced sheaf structures and the geometry of the simplicial complex is one of the directions for future research.

\paragraph{Cellular Sheaves of Modules}
As a generalization of cellular sheaves of ideals of a given ring \( R \) on a simplicial complex \( K \), one can also consider cellular sheaves of \( R \)-modules, that is, functors of the form \( \mathcal{F} : (K, \leq) \rightarrow \textup{\textsf{Mod}}_R \). However, inspired by weighted (persistent) homology~\cite{dawson1990homology,bell2019weighted,ren2018weighted,meng2019weighted}---particularly the mathematical framework presented in~\cite{dawson1990homology,meng2019weighted}, which has been applied in biomolecular data analysis~\cite{meng2019weighted,anand2020weighted}---we consider cellular sheaves of the form \( \mathcal{F} : (K, \leq)^{\operatorname{op}} \rightarrow \textup{\textsf{Mod}}_R \), where \( (K, \leq)^{\operatorname{op}} \) denotes the opposite category of \( (K, \leq) \). In this setting, we show that the weight assignments on the simplicial complex can be naturally interpreted via the sheaf cohomology of \( \mathcal{F} \).

\begin{remark}
In some literature, a functor of the form \( (K, \leq)^{\operatorname{op}} \rightarrow \textup{\textsf{Mod}}_R \) is typically referred to as a \textbf{cellular cosheaf} on the simplicial complex \( K \) (see, for example,~\cite{curry2014sheaves,curry2015topological}). However, to reduce notational overhead, we continue to use the term \textbf{cellular sheaf} for functors of this form throughout this paper, since \( (K, \leq)^{\operatorname{op}} \) remains a poset and our focus is on cellular sheaves defined over partially ordered sets.   
\end{remark}

\begin{def.}[\cite{ren2018weighted}]
Let \( (K, \leq) \) be a simplicial complex with vertex set \( V = \{ 1, 2, \ldots, n \} = K_{(0)} \). A \textbf{weight} on \( K \) with values in a ring \( R \) is a function \( w : K \rightarrow R \) satisfying the following condition: \( w(\sigma) \mid w(\tau) \) whenever \( \sigma \leq \tau \).
\end{def.}

Given a weight function \( w: K \rightarrow R \), a \textit{weighted boundary map} \( \partial^w_q: C_q(K; R) \rightarrow C_{q-1}(K; R) \) for \( q \in \mathbb{N} \) is defined as follows. Let \( C_q(K; R) := \bigoplus_{\sigma \in K_{(q)}} R\sigma \) denote the free \( R \)-module generated by all \( q \)-simplices. Then, for any \( q \)-simplex \( \sigma = [i_0, i_1, \ldots, i_q] \) with orientation \( i_0 < i_1 < \cdots < i_q \), the weighted boundary map is given by
\begin{equation}\label{Eq. weighted boundary map}
\partial_q^w(\sigma) = \sum_{k = 0}^q~(-1)^k \cdot \frac{w(\sigma)}{w(d_k(\sigma))} \cdot d_k(\sigma),
\end{equation}
where \( d_k(\sigma) := [i_0, \ldots, \widehat{i_k}, \ldots, i_q] \) denotes the \( k \)-th face of \( \sigma \). With the convention that \( \partial_0^w \) is the zero map, it holds that \( \partial_{q-1}^w \circ \partial_q^w = 0 \) for every \( q \in \mathbb{N} \) (see, e.g.,~\cite{dawson1990homology,ren2018weighted}). Using the incidence function from sheaf cohomology, the form of the weighted boundary map in Equation \eqref{Eq. weighted boundary map} can be reformulated by
\begin{equation}\label{Eq. weighted boundary map-2}
\partial_q^w(\sigma) = \sum_{k = 0}^q~ [d_k(\sigma):\sigma] \cdot \frac{w(\sigma)}{w(d_k(\sigma))} \cdot d_k(\sigma).
\end{equation}
Furthermore, in Equation~\eqref{Eq. weighted boundary map-2}, by viewing each ring element \( a \in R \) as an \( R \)-module homomorphism \( R \rightarrow R \), defined by scalar multiplication \( r \mapsto a \cdot r \), the weighted boundary map \( \partial_q^w \) can be interpreted as an \( R \)-linear map between free \( R \)-modules. Specifically, by defining a cellular sheaf as a functor \( \mathcal{F} : (K, \leq)^{\operatorname{op}} \rightarrow \textup{\textsf{Mod}}_R \), this framework can be interpreted in terms of sheaf coboundary maps and sheaf cohomology.

\begin{def.}\label{Definition: functor from Kop to ModR}
Let \( (K, \leq) \) be a simplicial complex with vertex set \( V = \{ 1, 2, \ldots, n \} = K_{(0)} \), and let \( w : K \rightarrow R \) be a weight on \( K \) with values in a ring \( R \). A functor \( \mathcal{F} : (K, \leq)^{\operatorname{op}} \rightarrow \textup{\textsf{Mod}}_R \) can be defined as follows:
\begin{itemize}
\item[\rm (a)] $\mathcal{F}_\sigma = R$ for every $\sigma \in K$.
\item[\rm (b)] $\mathcal{F}_{\tau,\sigma}: R \rightarrow R$ is defined as the $R$-module homomorphism $r \mapsto (w(\tau)/w(\sigma)) \cdot r$ for every pair $\sigma \leq \tau$ in $K$.
\end{itemize}
\end{def.}
We observe that $\mathcal{F}_{\rho,\sigma} = \mathcal{F}_{\tau,\sigma} \circ \mathcal{F}_{\rho,\tau}$ for every triple $\sigma \leq \tau \leq \rho$, as $(w(\rho)/w(\sigma)) = (w(\tau)/w(\sigma)) \cdot (w(\rho)/w(\tau))$. In other words,  \( \mathcal{F} : (K, \leq)^{\operatorname{op}} \rightarrow \textup{\textsf{Mod}}_R \) forms a functor. Then, Equation~\eqref{Eq. weighted boundary map-2} becomes
\begin{equation}\label{Eq. weighted boundary map-3}
\partial_q^w|_{R\sigma} = \sum_{k = 0}^q~ [d_k(\sigma):\sigma] \cdot \mathcal{F}_{\sigma, w(\sigma)},
\end{equation}
which is precisely the restriction of the $q$-th coboundary map \( \delta^q: C^q(K; \mathcal{F}) \rightarrow C^{q-1}(K; \mathcal{F}) \) to the subspace \( R_\sigma \), where the degree of the codomain is \( q-1 \) since \( \mathcal{F} \) is defined as a functor on the opposite category \( (K, \leq)^{\operatorname{op}} \). In other words, the weighted homology \( H_q(K; w) \) associated with a weight function \( w: (K, \leq) \rightarrow R \) coincides with the sheaf cohomology \( H^q(K; \mathcal{F}) \) of the cellular sheaf \( \mathcal{F} : (K, \leq)^{\operatorname{op}} \rightarrow \textup{\textsf{Mod}}_R \). We summarize this identification as the following theorem:

\begin{theorem}
Let \( (K, \leq) \) be a simplicial complex with vertex set \( V = \{ 1, 2, \ldots, n \} = K_{(0)} \), and let \( w : K \rightarrow R \) be a weight function with values in a ring \( R \). Let \( \mathcal{F} : (K, \leq)^{\operatorname{op}} \rightarrow \textup{\textsf{Mod}}_R \) be the cellular sheaf defined in Definition~\ref{Definition: functor from Kop to ModR}. Then,
\begin{equation*}
H_q(K; w) = H^q(K; \mathcal{F})    
\end{equation*}
for every \( q \in \mathbb{Z}_{\geq 0} \), where \( H_q(K; w) = \ker(\partial_q^w)/\operatorname{im}(\partial_{q+1}^w) \) is the \( q \)-th weighted homology.
\end{theorem}
 
In summary, cellular sheaves on simplicial complexes offer an alternative and flexible framework for encoding algebraic structures on the simplices of a complex. Through the sheaf-theoretic representation, various systematic constructions, such as direct sums, tensor products, and gluing—can be employed to generate new sheaf structures with enriched interactions across different dimensions. 

Moreover, leveraging the rich algebraic properties of rings, particularly graded rings such as the polynomial ring \( R = \mathbb{S}_n = \mathbb{F}[x_1, \ldots, x_n] \), the framework of cellular sheaves of \( R \)-modules has the potential to reveal deeper combinatorial and geometric features of the underlying simplicial complex, as well as data or signals supported on it. These capabilities provide a fertile ground for further theoretical development and practical applications in areas such as topological data analysis, combinatorial commutative algebra, and geometric learning. 
   


\section*{Acknowledgement}
HCS gratefully acknowledges the School of Physical and Mathematical Sciences (SPMS), Nanyang Technological University (NTU), Singapore, and thanks Dr. Kelin Xia (SPMS, NTU) for the inspiring research directions arising from their ongoing collaborations, which have motivated several of the questions explored in this paper. The initial motivation was also influenced by valuable discussions with Dr. Yu Wang and Dr. Ke Ye at the Academy of Mathematics and Systems Science, Chinese Academy of Sciences, China.

\bibliography{refs_2}
\bibliographystyle{abbrv}
\end{document}